\documentclass[a4paper,11pt]{amsart} 

\usepackage{amsmath,amssymb,amsxtra,amsthm}
\usepackage[margin=3cm]{geometry} 
\usepackage[hypertex]{hyperref}

\newtheorem{theorem}[equation]{Theorem}

\theoremstyle{remark}
\newtheorem{remark}[equation]{Remark}

\newcommand{\R}{\mathbb{R}} 
\newcommand{\N}{\mathbb{N}} 
\newcommand{\C}{\mathbb{C}} 

\newcommand{\cF}{\mathcal{F}} 
\newcommand{\PP}{\mathbb{P}} 
\newcommand{\hGamma}{\widehat{\Gamma}} 

\newcommand{\Res}{\operatorname{Res}} 

\newcommand{\iu}{\mathtt{i}} 

\begin{document} 

\title{Asymptotics of the Banana Feynman amplitudes at the 
large complex structure limit } 
\author{Hiroshi Iritani}  
\address{Department of Mathematics, Graduate School of 
Science, Kyoto University, Kitashirakawa-Oiwake-cho, 
Sakyo-ku, Kyoto, 6060-8502, Japan}
\email{iritani@math.kyoto-u.ac.jp}

\maketitle 
\begin{abstract} 
Recently B\"onisch-Fischbach-Klemm-Nega-Safari \cite{BFKNS} 
discovered, via numerical computation, that the leading 
asymptotics of the $l$-loop Banana Feynman amplitude 
at the large complex structure limit 
can be described by the Gamma class of 
a degree $(1,\dots,1)$ Fano hypersurface $F$ 
in $(\PP^1)^{l+1}$. 
We confirm this observation 
by using a Gamma-conjecture type result \cite{Iritani:periods} 
for $F$. 
\end{abstract}

\section{Introduction} 

The \emph{$l$-loop Banana Feynman amplitude} 
(see \cite[(8.1)-(8.2)]{Vanhove:exp}, \cite[(2.1)]{BFKNS}) 
is the integral 
\[
\cF(q,t) = \int_{(\R_{>0})^l} 
\frac{1}{t - \phi_q(y)} \frac{dy_1 \cdots dy_l}{y_1 \cdots y_l} 
\]
where $\phi_q$ is the 
Laurent polynomial 
\[
\phi_q(y) = (q_1 y_1 + \cdots+ q_l y_l + q_{l+1})
 (y_1^{-1}+\cdots+y_l^{-1} +1). 
\]
The parameters $q_i$, denoted by $\xi_i^2$ in \cite{BFKNS}, 
are the squares of the internal masses and $t$ is the square 
of the external momentum. 
When $t$ is a large positive number, 
the integrand has a pole along the hypersurface 
$(\phi_q(y) =t)$ 
and the integral diverges. We regularize the integral by means of 
analytic continuation: we know that the integral converges for $t<0$,  
and it can then be analytically continued to the complex plane (of $t$) 
minus the branch cut $[T, \infty)$, where 
$T := (\sum_{i=1}^{l+1} \sqrt{q_i} )^2 = 
\min\{\phi_q(y) : y\in (\R_{>0})^l\}$. 
 
The Feynman amplitude can be regarded as 
a relative period of the mixed Hodge structure of the pair 
$(\PP_\Delta \setminus M_{q,t}, \partial \PP_\Delta 
\setminus M_{q,t})$, where 
$\PP_\Delta$ is an $l$-dimensional toric variety  
such that $\phi_q^{-1}(t)\subset (\C^\times)^l$ 
is compactified to an anticanonical hypersurface $M_{q,t}\subset \PP_\Delta$ 
(intersecting every toric stratum properly) and 
$\partial \PP_\Delta = \PP_\Delta \setminus (\C^\times)^l$ 
is the toric boundary. 
As such, it satisfies \emph{inhomogeneous} Picard-Fuchs differential 
equations with respect to the parameters $q$ and $t$, 
which extend the Picard-Fuchs equations for 
$M_{q,t} = \overline{\phi_q^{-1}(t)}$. 
We refer the reader to \cite{Vanhove:exp, 
BKV:normal_fcn,BKV:local_sunset,KNS,BFKNS} 
and references therein 
for differential equations, 
Hodge-theoretic and arithmetic aspects 
of the Feynman amplitudes. 

In the present notes, we study the asymptotics of the Banana 
Feynman amplitude 
$\cF(q,t)$ near the large complex structure limit 
$t=\infty$ (or equivalently $q_1=\cdots = q_{l+1} = 0$) 
of $M_{q,t}$. 
\begin{theorem} 
\label{thm:main} 
Let $F$ be a degree $(1,\dots,1)$ Fano 
hypersurface in $(\PP^1)^{l+1}$ and 
let $p_1,\dots,p_{l+1}\in H^2(F)$ denote the hyperplane classes  
pulled back from $(\PP^1)^{l+1}$. 
For $q_1,\dots,q_{l+1}, t>0$, we have 
\[
\cF(q,t\mp \iu 0) 
\sim \frac{1}{t} \int_{F} e^{-p \log (q/t)} \cup 
\hGamma_{F} \cdot e^{\pm \pi \iu c_1(F)} \Gamma(1-c_1(F)) 
\quad \text{as $t\to \infty$} 
\]
where $p \log (q/t) = \sum_{i=1}^{l+1} p_i \log (q_i/t)$ and 
the sign depends on whether 
we perform the analytic continuation anti-clockwise 
or clockwise from a negative real $t$. 
\end{theorem} 

This follows from the power series expansion of $\cF(q,t)$ we give in 
Theorem \ref{thm:expansion} below. 
The class $\hGamma_F \in H^*(F)$ in the theorem 
is the \emph{Gamma class} \cite{Libgober,Iritani:Int,KKP} 
of the tangent bundle $TF$; it is explicitly given as 
\[
\hGamma_F = \frac{\Gamma(1+p_1)^2 \cdots \Gamma(1+p_{l+1})^2}
{\Gamma(1+p_1+\cdots+p_{l+1})} 
= \frac{e^{-2 \gamma c_1(F)}}{\Gamma(1+c_1(F))} 
\]
where $\Gamma(1+z)=\int_0^{\infty} 
e^{-t} t^z dt$ is the Euler $\Gamma$-function (when 
evaluating it at a cohomology class, we take its Taylor expansion) 
and $\gamma=0.577\cdots$ is the Euler constant. 
We also note that $c_1(F) = p_1+ \cdots + p_{l+1}$.

\begin{remark} 
Kerr \cite[Example 9.10]{Kerr} outlined another way to 
evaluate the asymptotics of the $l$-loop Banana Feynman integral. 
\end{remark} 

\begin{remark} 
Theorem \ref{thm:main} confirms the numerical computation 
by B\"onisch-Fischbach-Klemm-Nega-Safari \cite[\S 3]{BFKNS}. 
By taking the imaginary and the real parts of Theorem \ref{thm:main}, 
we get the following asymptotics as $t\to \infty$: 
\begin{align}
\nonumber 
\Im \cF(q,t-\iu 0) & \sim \frac{1}{t} \int_F  e^{-p \log (q/t)} \cup 
\frac{\prod_{j=1}^{l+1} \Gamma(1+p_j)^2}
{\Gamma(1+c_1(F))^2} \pi c_1(F) \\
\label{eq:ImF}
& = \frac{\pi}{t} \int_{W} e^{-p \log(q/t)} \cup \hGamma_W \\ 
\label{eq:ReF} 
\Re \cF(q,t-\iu 0) & \sim \frac{1}{t} \int_F 
e^{-p \log (q/t)} \cup 
\cos(\pi c_1(F)) 
\frac{\Gamma(1-c_1(F))}{\Gamma(1+c_1(F))} e^{-2 \gamma c_1(F)} 
\end{align} 
where $W \subset F$ is an anticanonical hypersurface, i.e.~the intersection 
of two degree $(1,\dots,1)$ hypersurfaces in $(\PP^1)^{l+1}$; 
this is a mirror of $M_{q,t}$. 
These asymptotics coincide\footnote{The factor $e^{-2\gamma c_1(F)}$ 
is missing in the second expression of \cite[(3.37)]{BFKNS}.} 
with \cite[(3.18); (3.36), (3.37)]{BFKNS}. 
\end{remark} 

\begin{remark} 
By the reflection principle, 
the imaginary part of $\cF(q,t-\iu 0)$ 
with $t>0$ can be understood 
as the difference $\frac{1}{2\iu}(\cF(q,t-\iu 0) - \cF(q,t+\iu 0))$ 
of two analytic continuations. 
This can then be equated with the residue integral 
\[
\pi \int_{\phi_q^{-1}(t) \cap (\R_{>0})^l} 
\Res\left( \frac{1}{t-\phi_q(y)} \frac{dy_1\cdots dy_l}{y_1\cdots y_l}
\right) 
\]
over the vanishing cycle $\phi_q^{-1}(t) \cap (\R_{>0})^l 
\subset M_{q,t}$. The Calabi-Yau Gamma conjecture 
\cite{Hosono,Iritani:periods} 
predicts that the asymptotics of such vanishing periods 
should be given by the Gamma class of the mirror partner $W$ 
of $M_{q,t}$, as in \eqref{eq:ImF}; in the case at hand 
this has been proved in \cite[Theorem 5.7]{Iritani:periods}. 
On the other hand, the asymptotics \eqref{eq:ReF} of the real part of $\cF(q,t)$ 
discovered in \cite{BFKNS} 
is slightly beyond the scope of the Calabi-Yau Gamma conjecture; 
it is related to (a degeneration of) the mixed Hodge structure 
(see also the recent work \cite{GKS}). 
In this paper, we will derive this from the \emph{Fano} 
Gamma conjecture \cite{Iritani:Int, KKP, Golyshev:deresonating, 
Iritani:periods, GGI}. 
\end{remark} 

\begin{remark} 
We can interpret $\hGamma_F \cdot \Gamma(1-c_1(F))$ as the 
Gamma class of the total space $K_F$ of the canonical bundle of $F$. 
See \cite{BKV:local_sunset} for the relation to local mirror symmetry. 
\end{remark} 

\section{Proof of the asymptotics} 
Theorem \ref{thm:main} follows immediately from the following result  
(compare \cite[Proposition 5.1]{Iritani:periods}). 
\begin{theorem} 
\label{thm:expansion} 
Let $q_1,\dots,q_{l+1}$ be positive real numbers. 
For $t \ll 0$, we have 
\[
\cF(q,t) = \frac{1}{t}\int_{F} I_W(q/(-t),-1) 
\cup \hGamma_{F} \cdot \Gamma(1-c_1(F)) 
\] 
where $I_W(q,z)$ is the cohomology-valued hypergeometric series 
\[
I_W(q,z) = e^{p\log q/z} \sum_{d=(d_1,\dots,d_{l+1}) 
\in \N^{l+1}} 
\frac{\prod_{k=1}^{d_1+\cdots+d_{l+1}} (p_1+\cdots+p_{l+1}+kz)^2}
{\prod_{i=1}^{l+1} \prod_{k=1}^{d_i} (p_i+ kz)^2} q^d 
\]
with $p \log q = \sum_{i=1}^{l+1} p_i \log q_i$ and 
$q^d = q_1^{d_1} \cdots q_{l+1}^{d_{l+1}}$. 
\end{theorem} 

\begin{remark} 
The hypergeometric series $I_W(q,z)$ is the Givental $I$-function 
\cite{Givental:toric} 
for the anticanonical hypersurface $W \subset F$, which is mirror to $M_{q,t}$. 
Here we regard it as a function taking values in $H^*(F)$, rather than 
in $H^*(W)$.  
\end{remark} 
A crucial observation \cite[p.41]{Vanhove:slide} 
is the fact that the Laurent polynomial $\phi_q(y)$ is a mirror 
of the Fano manifold $F$. 
The Givental mirror \cite[p.150, equation ($**$)]{Givental:toric} 
of the $(1,\dots,1)$-hypersurface 
$F\subset (\PP^1)^{l+1}$ is given by the oscillatory 
integral 
\[
\int_{C \subset \{u_1+\cdots + u_{l+1}=1\}} 
e^{-(\frac{q_1}{u_1} + \cdots + \frac{q_{l+1}}{u_{l+1}})} 
\frac{d \log u_1 \cdots d\log u_{l+1}}{d(u_1+\cdots +u_{l+1})}.  
\]
By the Przyjalkowsky change of variables \cite{Prz} 
\[
u_1 = \frac{y_1}{1+y_1+\cdots +y_l}, \quad 
\cdots \quad u_l = \frac{y_l}{1+y_1+\cdots +y_l}, \quad 
u_{l+1} = \frac{1}{1+y_1+\cdots + y_l}, 
\]
the above oscillatory integral can be rewritten as 
\[
\int_{C'} e^{-(1+y_1+\cdots+y_l) (\frac{q_1}{y_1} + 
\cdots + \frac{q_l}{y_l} + q_{l+1})} \frac{dy_1\cdots dy_l}{y_1\cdots y_l}. 
\]
The phase function equals the Laurent polynomial $-\phi_q(y)$ 
after the change of variables $y_i \to y_i^{-1}$. 
Therefore, the Gamma-conjecture type result 
\cite[Theorem 5.7]{Iritani:periods} implies that we have 
\begin{equation} 
\label{eq:gammathm} 
\int_{(\R_{>0})^l} e^{-\phi_q(y)} \frac{dy_1\cdots dy_l}{y_1\cdots y_l} 
= \int_{F} I_{F}(q,-1) \cup \hGamma_{F} 
\end{equation} 
for $q_1,\dots,q_{l+1}>0$, 
where $I_F$ is the Givental $I$-function \cite{Givental:toric} for $F$ 
\[
I_{F}(q,z) = e^{p \log q/z}
\sum_{d=(d_1,\dots,d_{l+1}) \in \N^{l+1}}
\frac{\prod_{k=1}^{d_1+\cdots+d_{l+1}} (p_1+\cdots+p_{l+1}+kz)}
{\prod_{i=1}^{l+1} \prod_{k=1}^{d_i} (p_i+ kz)^2} q^d.  
\] 
We substitute $r q_i$ for $q_i$ in the equation \eqref{eq:gammathm} 
and perform the Laplace transformation with respect to $r$. 
We find 
\begin{equation} 
\label{eq:Lap} 
\int_0^\infty e^{rt} dr \int_{(\R_{>0})^l} e^{-\phi_{rq}(y)} 
\frac{dy_1\cdots dy_l}{y_1\cdots y_l}  
= \int_0^\infty e^{rt} dr \int_{F} I_{F}(rq,-1) \cup \hGamma_{F} 
\end{equation} 
for $t<0$. 
A similar computation appeared in \cite[Section 5.1]{Iritani:periods}. 
Using $\phi_{rq}(y) = r \phi_q(y)$ and 
performing the integration in $r$ first\footnote
{This is legitimate, because the integrand is a positive continuous function.}, 
the left-hand side yields the Feynman amplitude 
\[ 
\int_{(\R_{>0})^l} 
\left( \int_0^\infty e^{(t-\phi_q(y)) r} dr\right) 
\frac{dy_1\cdots dy_l}{y_1\cdots y_l} 
= - \int_{(\R_{>0})^l} \frac{1}{t-\phi_q(y)} 
\frac{dy_1\cdots dy_l}{y_1\cdots y_l} 
= - \cF(q,t). 
\]
The right-hand side can be computed termwise, using 
\[
\int_0^\infty e^{rt} \prod_{i=1}^{l+1} 
(rq_i)^{d_i-p_i} dr = 
\Gamma(1+\textstyle\sum_{i=1}^{l+1}(d_i-p_i)) \dfrac{q^{d-p}}
{(-t)^{1+\sum_{i=1}^{l+1} (d_i-p_i)}}.  
\]
Note that the coefficient of $q^d$ in the series $I_F(q,-1)$ has the 
norm bounded by $C^{1+|d|}/|d|!$ for some $C>1$, 
where $|d| = d_1+\dots +d_{l+1}$. 
From this it follows that, for a sufficiently negative $t\ll 0$, 
we can interchange the sum over $d$ and the integral  
and that the right-hand side of \eqref{eq:Lap} converges; 
in particular the left-hand side also does. 
This proves Theorem \ref{thm:expansion}. 

\subsection*{Acknowledgements} 
I thank Albrecht Klemm for inviting me to think about 
the Banana Feynman amplitudes and Matt Kerr for 
helpful comments. This research is supported 
by JSPS Kakenhi Grant Number 
16H06335, 16H06337, 20K03582.

\providecommand{\arxiv}[1]{\href{http://arxiv.org/abs/#1}{arXiv:#1}}

\end{document}